\documentclass{amsproc}%
\usepackage{amssymb}
\usepackage{amsmath}
\usepackage{amsfonts}
\usepackage{graphicx}%
\providecommand{\U}[1]{\protect\rule{.1in}{.1in}}
\theoremstyle{plain}

\newtheorem{corollary}{Corollary}

\newtheorem{definition}{Definition}

\newtheorem{lemma}{Lemma}

\newtheorem{problem}{Problem}
\newtheorem{proposition}{Proposition}
\newtheorem{remark}{Remark}

\newtheorem{theorem}{Theorem}
\numberwithin{equation}{section}
\begin{document}
\title[Ideal Knots of Constant Curvature]{Local Structure of Ideal Shapes of Knots, II\ \\Constant Curvature Case}
\author{Oguz C. Durumeric}
\address{Deaprtment of Mathematics\\
University of Iowa\\
Iowa City, Iowa 52245}
\email{odurumer@math.uiowa.edu}
\date{February 2006}

\begin{abstract}
The thickness, $NIR(K)$ of a knot or link $K$ is defined to be the radius of
the largest solid tube one can put around the curve without any self
intersections, which is also known as the normal injectivity radius of $K$.
For $C^{1,1}$ curves $K$, $NIR(K)=min\{\frac{1}{2}DCSC(K),\frac{1}%
{sup\kappa(K)})\}$, where $\kappa(K)$ is the generalized curvature, and the
double critical self distance $DCSD(K)$ is the shortest length of the segments
perpendicular to $K$ at both end points. The knots and links in ideal shapes
(or tight knots or links) belong to the minima of \emph{ropelength =
length/thickness} within a fixed isotopy class. In this article, we prove that
$NIR(K)=\frac{1}{2}DCSC(K)$, for every relative minimum $K$ of ropelength in
$\mathbf{R}^{n}$ for certain dimensions $n$, including $n=3$.

\end{abstract}
\maketitle

\section{Introduction}

In this article, the local structure of $C^{1,1}$ relatively extremal knots
and links in $\mathbf{R}^{n}$ will be studied, particularly the extremal knots
and links with maximal constant generalized curvature$.$ The non-constant
curvature case was studied in our earlier article [D2]. The thickness or the
normal injectivity radius $NIR(K,\mathbf{R}^{n})$ of a knotted curve (or link)
is the radius of the largest tubular neighborhood around the curve without
intersections of the normal discs. Several different notations for thickness
appeared in the literature. $R(K)$ was used for thickness in [LSDR] and [BS].
[GM] showed that the thickness $\eta_{\ast}(K)$ was equal to the minimum
$\Delta(K)$ of $\rho_{G}$, the global radius of curvature for $C^{2}$ curves.
In [CKS], Cantarella-Kusner-Sullivan defined thickness $\tau(K)$ by the
infimum of the global radius of curvature and proved that it was the normal
injectivity radius for $C^{1,1}$ curves.

The ideal knots are the embeddings of $S^{1}$ into $\mathbf{R}^{3},$
maximizing $NIR(K,\mathbf{R}^{3})$ in a fixed isotopy (knot) class of fixed
length. More generally, a relatively extremal knot is a relative minimum of
the ropelength or isoembolic length, $\ell_{e}(K)=\frac{\ell(K)}%
{NIR(K,\mathbf{R}^{n})}$ in $C^{1}$ topology, where $\ell$ denotes the usual
length. The tight links and ideal knots belong to the relative minima of the
ropelength. The notion of ropelength has been defined and studied by several
authors, Litherland-Simon-Durumeric-Rawdon (called its reciprocal thickness in
[LSDR]), Gonzales and Maddocks [GM], Cantarella-Kusner-Sullivan [CKS] and
others. Cantarella-Kusner-Sullivan [CKS] defined ideal (thickest) knots as
\textquotedblleft tight\textquotedblright\ knots.

As J. Simon pointed out that there are physical examples (no proofs) of
relatively extremal unknots in $\mathbf{R}^{3}$, which are not circles, and
hence not ideal knots. One can construct similar physical examples for
composite knots. For dimensions $n\neq3,$ every $1$-dimensional knot is
trivial through an isotopy of curves \emph{of zero thickness}. At a strict
relative minimum $K$ of ropelength, one can not isotope the $n$-dimensional
solid tube of radius $NIR(K)$ around $K$ without increasing the length of $K.
$ Hence, one should not assume that all of the relative minima of ropelength
in $\mathbf{R}^{n}$ (for $n\neq3$) is the absolute minimum, that is a planar circle.

The thickness can be written in terms of \ the generalized curvature $\kappa$
and double critical self distance $DCSD(K)$ which is the shortest length of
the segments perpendicular to $K$ at both end points. \emph{Section 2} has the
formal definitions. Thickness Formula was discussed for $C^{2}-$knots in
$\mathbf{R}^{3}$ in [LSDR], and for $C^{1,1}$ knots in $\mathbf{R}^{3}$ by
Litherland in [L]. Also, [CKS, \emph{Lemma 1}] proved the Thickness Formula
below for $C^{1,1}$ knots and links in $\mathbf{R}^{3},$ since the geometric
and analytic curvatures are the same: $F_{g}=F_{k}=1/(\sup\kappa)$ by [D2,
\emph{Lemma 2}].

The notion of the global radius of curvature $\rho_{G}$ developed by Gonzales
and Maddocks for smooth curves in $\mathbf{R}^{3}$ defined by using circles
passing through 3 points of the curve in [GM] is another characterization of
$NIR(K,\mathbf{R}^{3}).$ This is still true for all continuous curves by [CKS,
\emph{Lemma 1}]. The construction of $\rho_{G}$ and rolling ball radius
$R_{O}$ for curves in $\mathbf{R}^{3}$ are different in nature due to 3-point
intersection condition versus 1-point of tangency and 1-point of intersection
condition. However, at the infimum they tend towards the same quantity:
$NIR(K,\mathbf{R}^{3}),$ [CKS].

$NIR(K,M)=R_{O}(K,M),$ a rolling ball/bead description of the injectivity
radius in $\mathbf{R}^{n}$, was known by Nabutowsky for hypersurfaces, by Buck
and Simon for $C^{2}$ curves, [BS], and by Cantarella, Kusner and Sullivan
[CKS, \emph{Lemma 1}]. Although the equality $NIR(K,M)=R_{O}(K,M)$ is
generalizable to all dimensions and to Riemannian manifolds [D1], the notion
of $\rho_{G}$ can not be used beyond the spaces of constant curvature.

In all of our results, the manifolds $K$\ are allowed to have several
components (unless stated otherwise). If $K$\ is one dimensional, we will use
$\gamma:D\rightarrow K$ for a parametrization of $K$ where $D$\ can be taken
as a finite disjoint union of intervals and circles, $S^{1}.$ Hence, all
closed curves are $C^{1}$ at the closing point.

\bigskip

\textbf{GENERAL\ THICKNESS\ FORMULA }[D1, \emph{Theorem 1}]

\emph{For every complete smooth Riemannian manifold }$M^{n}$\emph{\ and every
compact }$C^{1,1}$\emph{\ submanifold }$K^{k}$\emph{\ }$(\partial
K=\emptyset)$\emph{\ of }$M,$\emph{\ }
\[
NIR(K,M)=R_{O}(K,M)=\min\{F_{g}(K),\frac{1}{2}DCSD(K)\}.
\]

\textbf{THICKNESS\ FORMULA in }$\mathbf{R}^{n}$, [CKS, L, D1, D2]: \emph{For
every union of finitely many disjoint }$C^{1,1}$ \emph{simple closed curves
}$K$\emph{\ in }$\mathbf{R}^{n}$\emph{, one has}
\[
NIR(K,\mathbf{R}^{n})=R_{O}(K,\mathbf{R}^{n})=\min\{F_{k}(K),\frac{1}%
{2}DCSD(K)\}.
\]

\begin{remark}
Since all problems we discuss in this article involve bounded curvature in
$\mathbf{R}^{n},$ we rescale and take $\sup\kappa=\Lambda=1$ to simplify our
statements and proofs.
\end{remark}

The main question we address in [D2] and this article is \emph{Problem 1}
which is closely related to \emph{Problem 2, Markov-Dubins Problem:}

\begin{problem}
For which relatively extremal knots and links for the ropelength functional in
$\mathbf{R}^{n}$ does one have $NIR(K,\mathbf{R}^{n})=\frac{1}{2}DCSD(K)?$
\end{problem}

\begin{problem}
(Markov [M]-Dubins [Du]) Given $p,q,v,w$ in $\mathbf{R}^{n},$ with $\left\Vert
v\right\Vert =\left\Vert w\right\Vert =1.$ Classify all of the shortest curves
in $\mathcal{C}(p,q;v,w)$ which is the set of all curves $\gamma$ between the
points $p$ and $q$ in $\mathbf{R}^{n}$ with $\gamma^{\prime}(p)=v$,
$\gamma^{\prime}(q)=w$ and $\kappa\gamma\leq1=\Lambda$.
\end{problem}

\begin{definition}
A $C^{1,1}$ curve $\gamma:I=[a,b]\rightarrow\mathbf{R}^{n}\ $is called \ a
$CLC-$curve if there are $a\leq c\leq d\leq b$ such that (a) $\gamma([c,d])$
is a line segment of possibly zero length, and (b) each of $\gamma([a,c])$ and
$\gamma([d,b])$ is a planar circular arc of radius $1$ and of length in
$[0,2\pi)$. The curve $\gamma$ need not be planar. Similarly, one can define
$CCC$-curves by $C^{1}-$concatenation of 3 arcs of circles of curvature $1$,
where each arc has positive length and successive arcs have different centers.
\end{definition}

In an earlier article [D2], the author resolved \emph{Problem 1} in
$\mathbf{R}^{n},$ if the curve $K$ did not have constant curvature, see below.
The main tool in proving [\emph{Theorem 2} of D2] is the study of
$CLC-$curves. It tells us that the sections of a relatively extremal knot in
$\mathbf{R}^{n}$ with the minimal double critical points $I_{c}(K)$ removed
are $CLC-$curves or overwound, i.e. $\kappa\equiv\Lambda.$ This generalizes
one of the earliest results about the shape of ideal knots, that was obtained
by Gonzales and Maddocks [GM, p 4771]: \emph{A smooth ideal knot can be
partitioned into arcs of constant (maximal) global curvature and line
segments. }$I_{c}(K)$ denotes the minimal DCSD points on $K$ below.

\bigskip

[D2, \textbf{Theorem 2}] \ \emph{Let }$K$\emph{\ be a union of finitely many
disjoint }$C^{1,1}$ \emph{simple closed curves in }$\mathbf{R}^{n}$\emph{\ and
}$\gamma:D\rightarrow K\subset\mathbf{R}^{n}$\emph{\ be a parametrization. If
}$K$\emph{\ is a relative minimum of }$\ell_{e}$\emph{\ and }$\exists s_{0}\in
D,\kappa\gamma(s_{0})<\sup\kappa\gamma$\emph{, then both of the following
holds.}

\emph{(a) }$NIR(K,\mathbf{R}^{n})=R_{O}(K,\mathbf{R}^{n})=\frac{1}{2}%
DCSD(K)$\emph{.}

\emph{(b) If }$s_{0}\notin I_{c}(K)$\emph{, then there exists }$a,b$%
\emph{\ such that }$s_{0}\in\lbrack a,b],$\emph{\ }$\gamma([a,b])$\emph{\ is a
CLC(}$\sup\kappa\gamma$\emph{)-curve where the line segment has positive
length and contains }$\gamma(s_{0})$\emph{, and each circular part has at most
}$\pi$\emph{\ radians angle ending at a point of }$I_{c}(K).$

[D2, \textbf{Corollary 2}] \emph{Let }$K$\emph{\ be a union of finitely many
disjoint }$C^{1,1}$\emph{\ simple closed curves in }$\mathbf{R}^{n}%
.$\emph{\ If }$K$\emph{\ is a relative minimum for }$\ell_{e}$\emph{\ and
curvature of }$K$\emph{\ is not identically constant }$R_{O}(K)^{-1}$\emph{,
then the thickness of }$K$\emph{\ is }$\frac{1}{2}DCSD(K)$\emph{.
Equivalently, if there exists a relative minimum }$K$\emph{\ for }$\ell_{e}%
$\emph{\ such that }$\frac{1}{2}DCSD(K)>R_{O}(K)=F_{k}(K),$\emph{\ then }%
$K$\emph{\ must have constant generalized curvature }$F_{k}(K)^{-1}.$

\bigskip

Markov [M], Dubins [Du] and Reeds and Shepp [ReSh] studied the $2-$dimensional
cases for \emph{Problem 2}. In dimension $3$, the following results of H.
Sussmann obtain the possible types solutions for this problem. A helicoidal
arc is a smooth curve in $\mathbf{R}^{3}$ with constant curvature $1$ and
positive torsion $\tau$ satisfying the differential equation $\tau
^{\prime\prime}=1.5\tau^{\prime}\tau^{-1}-2\tau^{3}+2\tau-\zeta\tau\left\vert
\tau\right\vert ^{1/2}$ for some nonnegative constant $\zeta.$

\bigskip

\textbf{THEOREM. }\emph{(Sussmann [S]) }

\emph{1. For the Markov-Dubins problem in dimension three, every minimizer is
either (a) a helicoidal arc or (b) a concatenation of three pieces each of
which is a circle or a straight line. For a minimizer of the form CCC, the
middle arc has length }$\geq\pi$ \emph{and }$<2\pi$.

\emph{2. Every helicoidal arc corresponding to a value of }$\zeta$ \emph{such
that }$\zeta>0$\emph{\ is local strict minimizer.}

\bigskip

Sussmann further proves that \emph{CSC-conjecture} (every minimizer is either
\emph{CCC} or \emph{CLC, }[ReSh]) is false in $\mathbf{R}^{3}$ [S,
\emph{Propositions 2.1} and \emph{2.2}]. In [S], the details of the steps of
the proof of \emph{Theorem 1} of Sussmann are provided, but there are only few
remarks about proof of its \emph{Theorem 2}.

\bigskip

The main result of this article is the following \emph{Theorem 1} below which
shows the nonexistence of a relative minimum $K$ for $\ell_{e}$ with $\frac
{1}{2}DCSD(K)>F_{k}(K)$, in all cases in certain dimensions, including
constant curvature cases. Only \emph{Corollary 1} invokes [S, \emph{Theorem
1}]. The remaining results of this article including \emph{Theorem 1} and [D2]
are independent of Sussmann [S].

\begin{theorem}
Let $n$\ be a dimension such that

(i) every minimizer for the Markov-Dubins problem in $\mathbf{R}^{n}$ is
either a smooth curve with curvature 1 and positive torsion, or a $C^{1}%
-$concatenation of finitely many circular arcs of curvature $1$ and a line
segment, and

(ii) every $CCC-$curve with the middle arc of length $<\pi$\ is not a minimizer.

Then, $NIR(K,\mathbf{R}^{n})=\frac{1}{2}DCSD(K)$ for every relative minimum
$K$ of $\ell_{e}$ where $K$ is a union of finitely many disjoint $C^{1,1}$
simple closed curves in\textbf{\ }$\mathbf{R}^{n}.$
\end{theorem}

\begin{corollary}
Let $K$ be a union of finitely many disjoint $C^{1,1}$ simple closed curves in
$\mathbf{R}^{3}.$ If $K$ is a relative minimum of $\ell_{e}$, then
$NIR(K,\mathbf{R}^{3})=\frac{1}{2}DCSD(K).$
\end{corollary}

\section{Definitions and Notation}

We assume all parametrizations $\gamma:D\rightarrow K\subset\mathbf{R}^{n}$
are one-to-one and $\left\Vert \gamma^{\prime}\right\Vert \neq0$.

\begin{definition}
$\exp_{p}^{N}v=p+v:NK\rightarrow$ $\mathbf{R}^{n}$ is the normal exponential
map of $K$ in $\mathbf{R}^{n}$. The thickness of $K$ in $\mathbf{R}^{n}$ or
the normal injectivity radius of $\exp^{N}$ is
\[
NIR(K,\mathbf{R}^{n})=\sup(\left\{  0\right\}  \cup\{r>0:\exp^{N}:\{v\in
NK:\left\Vert v\right\Vert <r\}\rightarrow M\text{ is one-to-one}\}).
\]
Equivalently, if $\gamma(s)$ parametrizes $K,$ then

$r>NIR(K,\mathbf{R}^{n})\Leftrightarrow\left(
\begin{array}
[c]{c}%
\exists\gamma(s),\gamma(t),q\in\mathbf{R}^{n},\\
\gamma(s)\neq\gamma(t),\left\Vert \gamma(s)-q\right\Vert <r,\left\Vert
\gamma(t)-q\right\Vert <r,\text{ and }\\
(\gamma(s)-q)\cdot\gamma^{\prime}(s)=(\gamma(t)-q)\cdot\gamma^{\prime}(t)=0
\end{array}
\right)  .$
\end{definition}

\begin{definition}
For $\gamma:I\rightarrow\mathbf{R}^{n}$, define:

Dilation: $dil^{\alpha}\gamma^{\prime}(s,t)=\frac{\measuredangle
(\gamma^{\prime}(s),\gamma^{\prime}(t))}{\ell(\gamma([s,t]))}$ for $s\neq t$

(Generalized) Curvature: $\kappa\gamma(s)=\underset{t\neq u\text{ and
}t,u\rightarrow s}{\lim\sup\text{ }}dil^{\alpha}\gamma^{\prime}(t,u)$

Analytic focal distance: $F_{k}(\gamma)=(\underset{}{\sup_{I}}\kappa
\gamma(s))^{-1}.$

If $K$\textbf{\ }is a union of finitely many disjoint $C^{1,1}$ curves
$\gamma_{(i)}$ in\textbf{\ }$\mathbf{R}^{n},$ then $F_{k}(K)=\min_{i}%
F_{k}(\gamma_{(i)}).$
\end{definition}

\begin{definition}
Let $K$ be a finite union of disjoint $C^{1}$ curves in $\mathbf{R}^{n}$. For
any $v\in UT\mathbf{R}_{p}^{n}$ and any $r>0,$ define

(a) $O_{p}(v,r)=\{x\in\mathbf{R}^{n}:\exists w\in\mathbf{R}^{n},v\cdot
w=0,\left\|  w\right\|  =1,\left\|  x-p-rw\right\|  <r\}$

(b) $O_{p}^{c}(v,r)=\mathbf{R}^{n}-O_{p}(v,r)$

(c) $O_{p}(r;K)=O_{p}(v,r)$ where $v\in UTK_{p}$

(d) $O(r;K)=%
{\displaystyle\bigcup\limits_{p\in K}}
O_{p}(r;K)$

In all of the above, $r$ may be omitted when $r=1.$ $K$ will be omitted unless
there is an ambiguity.
\end{definition}

\begin{definition}
Let $K$ be a finite union of $C^{1}$ curves in $\mathbf{R}^{n}.$ Define

(a) The ball radius of $K$ in $\mathbf{R}^{n}$ to be $R_{O}(K,\mathbf{R}%
^{n})=\inf\{r>0:O(r;K)\cap K\neq\emptyset\}$

(b) The pointwise geometric focal distance $F_{g}(p)=\inf\{r>0:p\in
\overline{O_{p}(r;K)\cap K}\}$ for any $p\in K,$ and the geometric focal
distance $F_{g}(K)=\inf_{p\in K}F_{g}(p).$
\end{definition}

\begin{definition}
A pair of distinct points $p$ and $q$ in $K$ are called a double critical pair
for $K,$ if the line segment $\overline{pq}$ is normal to $K$ at both $p $ and
$q.$ The double critical self distance is
\[
DCSD(K)=\inf\{\left\Vert p-q\right\Vert :\{p,q\}\text{ is a double critical
pair for }K\}
\]
A double critical pair $\{p,q\}$ is called minimal if $DCSD(K)=\left\Vert
p-q\right\Vert .$
\end{definition}

\section{Review of Some Basic Tools from [D2]}

$K$\ denotes a union of finitely many disjoint $C^{1,1}$\ simple closed curves
in $\mathbf{R}^{n}$ and $\gamma:D\rightarrow K\subset\mathbf{R}^{n}$\ denotes
a one-to-one non-singular parametrization,\ where $D=%
{\textstyle\bigcup\nolimits_{i=1}^{k}}
S_{(i)}^{1}$, a union of $k$\ copies disjoint circles, unless stated
otherwise. When $\left\Vert \gamma^{\prime}\right\Vert \equiv1$ is assumed,
$S_{(i)}^{1}$ are taken with the appropriate radius and length. A knot or link
class $[\theta]$ is a free $C^{1}$ (ambient) isotopy class of embeddings of
$\gamma:D\rightarrow\mathbf{R}^{n}$ with a fixed number of components. Since
all of our proofs involve local perturbations of only one component at a time,
we will work with $\gamma_{(i)}:\mathbf{S}_{(i)}^{1}\rightarrow\mathbf{R}^{n}$
and we will omit the lower index $_{(i)}$ to simplify the notation wherever it
is possible. We will identify $\mathbf{S}^{1}\cong\mathbf{R}/L\mathbf{Z}$, for
$L>0,$ and use interval notation to describe connected proper subsets of
$\mathbf{R}/L\mathbf{Z.}$ In other words, $\gamma_{(i)}(t+L)=\gamma_{(i)}(t) $
and $\gamma_{(i)}^{\prime}(t+L)=\gamma_{(i)}^{\prime}(t)$, $\forall
t\in\mathbf{R}$ with $\left\Vert \gamma_{(i)}^{\prime}\right\Vert \neq0$ and
$\gamma_{(i)}$ is one-to-one on $[0,L).$ See [D2] for proofs of the following
propositions that will be used in this article.

\begin{lemma}
Let $\gamma:D\rightarrow K\subset\mathbf{R}^{n}$ be a $C^{1}$ knot or link.

(a) If $DCSD(K)>0,$ then there exists a critical pair $\{p_{0},q_{0}\}$ such that

$DCSD(K)=\left\Vert p_{0}-q_{0}\right\Vert .$

(b) If $\sup\kappa\gamma<\infty,$ i.e. $\gamma$ is $C^{1,1}$, then
$DCSD(K)>0.$
\end{lemma}

\begin{proposition}
Let $\left\{  \gamma_{m}\right\}  _{m=1}^{\infty}:D\rightarrow\mathbf{R}^{n} $
be a sequence uniformly converging to $\gamma$ in $C^{1}$ sense, i.e.
$(\gamma_{m}(s),\gamma_{m}^{\prime}(s))\rightarrow$ $(\gamma(s),\gamma
^{\prime}(s))$ uniformly on $D.$ Let $K_{m}=\gamma_{m}(D)$ for $m\geq1$ and
$K=\gamma(D).$

(a) ([CKS, Lemma 3] and [L]) If $R_{O}(K_{m})\geq r$ for sufficiently large
$m$, then $R_{O}(K)\geq r$. Consequently, $\lim\sup_{m}R_{O}(K_{m})\leq
R_{O}(K).$

(b) If $\lim\inf_{m}DCSD(K_{m})>0,$ then $\lim\inf_{m}DCSD(K_{m})\geq
DCSD(K).$
\end{proposition}

\begin{definition}
For $\gamma:D\rightarrow K\subset\mathbf{R}^{n},$ define

(a) $I_{c}=\{x\in D:\exists y\in D$ such that $\left\Vert \gamma
(x)-\gamma(y)\right\Vert =DCSD(K)$ and

$\qquad\left(  \gamma(x)-\gamma(y)\right)  \cdot\gamma^{\prime}(x)=\left(
\gamma(x)-\gamma(y)\right)  \cdot\gamma^{\prime}(y)=0\}$ and $K_{c}=\gamma
_{c}=\gamma(I_{c})$

(b) $I_{z}=\{x\in D:\kappa\gamma(x)=0\}$ and $K_{z}=\gamma_{z}=\gamma(I_{z})$

(c) $I_{mx}=\{x\in D:\kappa\gamma(x)=1/R_{O}(K)\}$ and $K_{mx}=\gamma
_{mx}=\gamma(I_{mx})$

(d) $I_{b}=\{x\in D:0<\kappa\gamma(x)<1/R_{O}(K)\}$ and $K_{b}=\gamma
_{b}=\gamma(I_{b})$
\end{definition}

\begin{proposition}
([CKS, Theorem 7], [GL], [GMSM]) For any knot/link class $[\theta]$ in
$\mathbf{R}^{n}$, $\exists\gamma_{0}\in\lbrack\theta]$ such that

(a) $\forall\gamma\in\lbrack\theta],$ $0<\ell_{e}(\gamma_{0})\leq\ell
_{e}(\gamma),$ and hence

(b) $\forall\gamma\in\lbrack\theta],$ $\left(  \ell(\gamma_{0})=\ell
(\gamma)\Longrightarrow R_{O}(\gamma_{0})\geq R_{O}(\gamma)\right)  .$
\end{proposition}

\begin{proposition}
Let $\left\{  \gamma_{m}\right\}  _{m=1}^{\infty}:D\rightarrow\mathbf{R}^{n} $
be a sequence uniformly converging to $\gamma$ in $C^{1}$ sense, $K=\gamma(D)$
and $K_{m}=\gamma_{m}(D),$ such that $\exists C<\infty,\forall m,$ $\sup
\kappa\gamma_{m}\leq C.$

(a) Let $A\subset D$ be a given compact set with $\{s\in D:\gamma_{m}%
(s)\neq\gamma(s)\}\subset A,\forall m.$ If $A\cap I_{c}=\emptyset,$ then
$\exists m_{1}$ such that $\forall m\geq m_{1},DCSD(K_{m})\geq DCSD(K).$

(b) If $F_{k}(K)<\frac{1}{2}DCSD(K)$ and $F_{k}(K_{m})\geq F_{k}(K),\forall
m,$ then $\exists m_{1}$ such that $\forall m\geq m_{1},R_{O}(K_{m})\geq
R_{O}(K).$
\end{proposition}

\begin{proposition}
(Also see [GM, p4771] for another version for smooth ideal knots.) Let $K$ be
a $C^{1,1}$ relatively minimal knot or link for the ropelength $\ell_{e}$.

(a) If $DCSD(K)=2R_{O}(K),$ then $K-(K_{c}\cup K_{mx})$ is a countable union
of open ended line segments, and hence $I_{b}\subset I_{c}.$

(b) If $DCSD(K)>2R_{O}(K)$, then $K-K_{mx}$ is a countable union of open ended
line segments, (in fact $\emptyset$ by [\emph{Theorem 2} of D2])
\end{proposition}

\section{Proof of Theorem 1}

In dimension 3, the following coincides with the standard definitions except
the sign of the torsion. For a $C^{3}$ curve $\gamma:I\rightarrow
\mathbf{R}^{n}$, $n\geq3,$ parametrized by arclength, i.e. $\left\Vert
\gamma^{\prime}(t)\right\Vert =1$, define

(a) $\mathbf{T=}\gamma^{\prime}(t)$,

(b) $\kappa=\left\Vert \mathbf{T}^{\prime}\right\Vert ,$ and if $\kappa>0, $
define $\mathbf{N=}\frac{1}{\kappa}\mathbf{T}^{\prime},$

(c) $\tau=\left\Vert \mathbf{N}^{\prime}+\kappa\mathbf{T}\right\Vert \geq0,$
and if $\tau>0,$ define $\mathbf{B=}\frac{1}{\tau}\left(  \mathbf{N}^{\prime
}+\kappa\mathbf{T}\right)  .$

This definition yields an orthonormal set $\{\mathbf{T,N,B}\}$ along $\gamma,$
if both $\kappa,\tau>0.$

\begin{proposition}
Let $K$ be a union of finitely many disjoint $C^{1,1}$ simple closed curves in
$\mathbf{R}^{n}$ and $K$ be a relative minimum of $\ell_{e}.$ If $K$ has a
component $K_{0}$ which is a $C^{4}$ simple closed curve of positive torsion
$\tau>0$ everywhere, then $NIR(K,\mathbf{R}^{n})=\frac{1}{2}DCSD(K)$.
\end{proposition}

\begin{proof}
(a) First prove the statement for a connected $K.$ Let $\gamma:\mathbf{S}%
^{1}\rightarrow K\subset\mathbf{R}^{n}$ parametrize $K$. Proposition holds if
curvature of $\gamma$ is not identically $NIR(K,\mathbf{R}^{n})^{-1}$, by [D2,
\emph{Corollary 2}].

By rescaling, assume that $\kappa\gamma\equiv NIR(K,\mathbf{R}^{n})^{-1}=1. $%
\begin{align*}
\mathbf{N}^{\prime}  &  =-\kappa\mathbf{T+}\tau\mathbf{B=}-\mathbf{T+}%
\tau\mathbf{B}\\
\mathbf{T}\cdot\mathbf{N}^{\prime}  &  =-\mathbf{T}^{\prime}\cdot
\mathbf{N}=-\kappa=-1\\
\mathbf{N}\cdot\mathbf{N}^{\prime\prime}  &  =\frac{1}{2}(\mathbf{N}%
\cdot\mathbf{N)}^{\prime\prime}-\mathbf{N}^{\prime}\cdot\mathbf{N}^{\prime
}\mathbf{=-}\left\Vert -\mathbf{T+}\tau\mathbf{B}\right\Vert ^{2}=-(1+\tau
^{2})
\end{align*}
Consider the variation $\gamma_{\varepsilon}(t)=\gamma(t)+\varepsilon
\mathbf{N}(t)$ and $\Gamma_{\varepsilon}(t)=\frac{L}{\ell(\gamma_{\varepsilon
})}\gamma_{\varepsilon}(t).$ By \emph{Lemma 2} below,
\begin{align*}
\left.  \frac{d}{d\varepsilon}\kappa\Gamma_{\varepsilon}\right\vert
_{\varepsilon=0}(t)  &  =\gamma^{\prime\prime}(t)\cdot\mathbf{N}^{\prime
\prime}(t)-2\gamma^{\prime}(t)\cdot\mathbf{N}^{\prime}(t)+\frac{1}{L}%
{\displaystyle\int\nolimits_{0}^{L}}
\gamma^{\prime}(u)\cdot\mathbf{N}^{\prime}(u)du\\
&  =\mathbf{N}\cdot\mathbf{N}^{\prime\prime}(t)-2\mathbf{T}\cdot
\mathbf{N}^{\prime}(t)+\frac{1}{L}%
{\displaystyle\int\nolimits_{0}^{L}}
\mathbf{T}\cdot\mathbf{N}^{\prime}(u)du\\
&  =-(1+\tau^{2})+2+\frac{1}{L}%
{\displaystyle\int\nolimits_{0}^{L}}
-1du\\
\left.  \frac{d}{d\varepsilon}\kappa\Gamma_{\varepsilon}\right\vert
_{\varepsilon=0}(t)  &  =-\tau^{2}(t)<0
\end{align*}
Since $K$ is compact and $\kappa\Gamma_{\varepsilon}$ is a $C^{2}$ function of
$t$ and $\varepsilon,$ there exists $\varepsilon_{0}>0$ such that
$\forall\varepsilon\in(0,\varepsilon_{0}),\forall t\in\mathbf{S}^{1}%
,(\kappa\Gamma_{\varepsilon}(t)<1)$. Hence, $\forall\varepsilon\in
(0,\varepsilon_{0})$, $\max\kappa\Gamma_{\varepsilon}<1$ and $F_{k}%
(\Gamma_{\varepsilon})>1=$ $F_{k}(K).$ Obviously, $\Gamma_{1/m}\rightarrow
\gamma$ in $C^{1}$ sense, as $m\rightarrow\infty.$

Suppose that $\frac{1}{2}DCSD(K)>F_{k}(K).$ Let $K_{m}=\Gamma_{1/m}%
(\mathbf{S}^{1}).$ \emph{By Proposition 1(b)},
\[
\lim\inf_{m}\frac{1}{2}DCSD(K_{m})\geq\frac{1}{2}DCSD(K)>F_{k}(K)=1
\]

For sufficiently large $m,$%
\begin{align*}
\frac{1}{2}DCSD(K_{m})  &  >1\\
F_{k}(K_{m})  &  >1\\
NIR(K_{m},\mathbf{R}^{n})  &  >1=F_{k}(K)=NIR(K,\mathbf{R}^{n})\text{
(Thickness Formula)}\\
\ell_{e}(K_{m})  &  <\ell_{e}(K),\text{ since }\ell(K_{m})=L=\ell(K)
\end{align*}
This contradicts to the fact that $K$ is relative minimum of $\ell_{e}.$

Consequently, $\frac{1}{2}DCSD(K)\leq F_{k}(K),$ that is $NIR(K,\mathbf{R}%
^{n})=\frac{1}{2}DCSD(K).$

(b) To prove this statement for $K$ with two or more components, let $K_{0}$
be a $C^{4}$ simple closed curve component of positive torsion $\tau>0$
everywhere. By [\emph{Corollary 2} of D2], the only remaining case is when all
of $K$ has constant curvature $1$. Let $\gamma:\mathbf{S}^{1}\rightarrow
K_{0}\subset\mathbf{R}^{n},$ and take the variation and rescaling
$\Gamma_{\varepsilon}(t)=\frac{\ell(K_{0})}{\ell(\gamma_{\varepsilon})}%
\gamma_{\varepsilon}(t)$ \textbf{only }along $K_{0}$ and leave the other
components invariant. Let $K_{m}^{\ast}=\Gamma_{1/m}(\mathbf{S}^{1}%
)\cup(K-K_{0}).$ Obviously, $K_{m}^{\ast}\rightarrow K$ in $C^{1}$ sense.
Repeat the same proof as in (a) until and including \textquotedblleft Suppose
that $\frac{1}{2}DCSD(K)>F_{k}(K)$\textquotedblright. By \emph{Proposition
1(b)},
\[
\lim\inf_{m}\frac{1}{2}DCSD(K_{m}^{\ast})\geq\frac{1}{2}DCSD(K)>F_{k}(K)=1
\]
For sufficiently large $m,$%
\begin{align*}
\frac{1}{2}DCSD(K_{m}^{\ast})  &  >1\\
F_{k}(K_{m}^{\ast})  &  =F_{k}(K-K_{0})=1<F_{k}(\Gamma_{1/m})\\
NIR(K_{m}^{\ast},\mathbf{R}^{n})  &  =1=F_{k}(K)=NIR(K,\mathbf{R}^{n})\text{
(Thickness Formula)}\\
\ell_{e}(K_{m}^{\ast})  &  =\ell_{e}(K),\text{ since }\ell(K_{m}^{\ast
})=L=\ell(K_{0})
\end{align*}
Since $K$ is a relative minimum of $\ell_{e},$ $K_{m}^{\ast}$ is a relative
minimum of $\ell_{e},$ for sufficiently large $m.$ $K_{m}^{\ast}$ does not
have constant maximal curvature $1$ everywhere, since $\max\kappa\Gamma
_{1/m}<1.$ Hence, by [D2, \emph{Corollary 2}]\emph{,} $1=NIR(K_{m}^{\ast
},\mathbf{R}^{n})=\frac{1}{2}DCSD(K_{m}^{\ast})$ for sufficiently large $m,$
which contradicts above. Consequently, $\frac{1}{2}DCSD(K)\leq F_{k}(K),$ that
is $NIR(K,\mathbf{R}^{n})=\frac{1}{2}DCSD(K).$
\end{proof}

\begin{lemma}
Let $\gamma:[0,L]\rightarrow\mathbf{R}^{n}$ be be a $C^{2}-$curve parametrized
by arclength with constant curvature $1$: $\left\Vert \gamma^{\prime
}(t)\right\Vert =\left\Vert \gamma^{\prime\prime}(t)\right\Vert =1$ and let
$V:[0,L]\rightarrow\mathbf{R}^{n}$ be a $C^{2}$ vector field along $\gamma$
normal to $\gamma$, i.e. $V(t)\cdot\gamma^{\prime}(t)=0.$ Define
$\gamma_{\varepsilon}(t)=\gamma(t)+\varepsilon V(t)$ and $\Gamma_{\varepsilon
}(t)=\frac{L}{\ell(\gamma_{\varepsilon})}\gamma_{\varepsilon}(t),$ where
$L=\ell(\gamma).$ Then

(a)$\left.  \frac{d}{d\varepsilon}\kappa\gamma_{\varepsilon}(t)\right\vert
_{\varepsilon=0}=\gamma^{\prime\prime}(t)\cdot V^{\prime\prime}(t)-2\gamma
^{\prime}(t)\cdot V^{\prime}(t)$ and

(b) $\left.  \frac{d}{d\varepsilon}\kappa\Gamma_{\varepsilon}(t)\right\vert
_{\varepsilon=0}=\gamma^{\prime\prime}(t)\cdot V^{\prime\prime}(t)-2\gamma
^{\prime}(t)\cdot V^{\prime}(t)+\frac{1}{L}%
{\displaystyle\int\nolimits_{0}^{L}}
\gamma^{\prime}(u)\cdot V^{\prime}(u)du$
\end{lemma}

\begin{proof}
We include this elementary computation for the sake of completeness.

(a) Recall that:
\[
\kappa\alpha=\left\Vert \alpha^{\prime\prime}(\alpha^{\prime}\cdot
\alpha^{\prime})-\alpha^{\prime}(\alpha^{\prime\prime}\cdot\alpha^{\prime
})\right\Vert \left\Vert \alpha^{\prime}\right\Vert ^{-4}%
\]%
\[
\left.  \frac{d}{d\varepsilon}\left\Vert v_{0}+\varepsilon v_{1}%
+\varepsilon^{2}v_{2}\right\Vert ^{k}\right\vert _{\varepsilon=0}=k(v_{0}\cdot
v_{1})\left\Vert v_{0}\right\Vert ^{k-2}%
\]%
\begin{align*}
\gamma^{\prime\prime}\cdot\gamma^{\prime}  &  =0\\
\left\Vert \gamma^{\prime}\right\Vert  &  =\left\Vert \gamma^{\prime\prime
}\right\Vert =\kappa\gamma=1\\
\gamma_{\varepsilon}^{\prime}(t)  &  =\gamma^{\prime}(t)+\varepsilon
V^{\prime}(t)\\
\gamma_{\varepsilon}^{\prime\prime}(t)  &  =\gamma^{\prime\prime
}(t)+\varepsilon V^{\prime\prime}(t)
\end{align*}%
\begin{align*}
\gamma_{\varepsilon}^{\prime\prime}(\gamma_{\varepsilon}^{\prime}\cdot
\gamma_{\varepsilon}^{\prime})-\gamma_{\varepsilon}^{\prime}(\gamma
_{\varepsilon}^{\prime\prime}\cdot\gamma_{\varepsilon}^{\prime})  &
=\gamma^{\prime\prime}+\varepsilon\left[  V^{\prime\prime}+2\gamma
^{\prime\prime}(V^{\prime}\cdot\gamma^{\prime})-\gamma^{\prime}(V^{\prime
\prime}\cdot\gamma^{\prime})-\gamma^{\prime}(\gamma^{\prime\prime}\cdot
V^{\prime})\right]  +o(\varepsilon^{2})\\
&  :=w_{0}+\varepsilon w_{1}+o(\varepsilon^{2})
\end{align*}%
\[
w_{0}\cdot w_{1}=\gamma^{\prime\prime}\cdot\left[  V^{\prime\prime}%
+2\gamma^{\prime\prime}(V^{\prime}\cdot\gamma^{\prime})-\gamma^{\prime
}(V^{\prime\prime}\cdot\gamma^{\prime})-\gamma^{\prime}(\gamma^{\prime\prime
}\cdot V^{\prime})\right]  =\gamma^{\prime\prime}\cdot V^{\prime\prime
}+2(V^{\prime}\cdot\gamma^{\prime})
\]%
\begin{align*}
\left.  \frac{d}{d\varepsilon}\kappa\gamma_{\varepsilon}\right\vert
_{\varepsilon=0}  &  =\left.  \frac{d}{d\varepsilon}\left\Vert \gamma
_{\varepsilon}^{\prime\prime}(\gamma_{\varepsilon}^{\prime}\cdot
\gamma_{\varepsilon}^{\prime})-\gamma_{\varepsilon}^{\prime}(\gamma
_{\varepsilon}^{\prime\prime}\cdot\gamma_{\varepsilon}^{\prime})\right\Vert
\left\Vert \gamma_{\varepsilon}^{\prime}\right\Vert ^{-4}\right\vert
_{\varepsilon=0}\\
&  =\left(  \gamma^{\prime\prime}\cdot V^{\prime\prime}+2(V^{\prime}%
\cdot\gamma^{\prime})\right)  \left\Vert \gamma^{\prime\prime}\right\Vert
^{-1}\left\Vert \gamma^{\prime}\right\Vert ^{-4}-4\left\Vert \gamma
^{\prime\prime}\right\Vert \left\Vert \gamma^{\prime}\right\Vert
^{-6}(V^{\prime}\cdot\gamma^{\prime})\\
&  =\gamma^{\prime\prime}\cdot V^{\prime\prime}+2V^{\prime}\cdot\gamma
^{\prime}-4V^{\prime}\cdot\gamma^{\prime}=\gamma^{\prime\prime}\cdot
V^{\prime\prime}-2V^{\prime}\cdot\gamma^{\prime}%
\end{align*}
(b)
\[
\Gamma_{\varepsilon}(t)=\frac{L}{\ell(\gamma_{\varepsilon})}\gamma
_{\varepsilon}(t)\text{ hence }\kappa\Gamma_{\varepsilon}(t)=\frac{\ell
(\gamma_{\varepsilon})}{L}\kappa\gamma_{\varepsilon}(t)
\]
By the classical First Variation Formula, [CE]:%
\begin{align*}
\left.  \frac{d}{d\varepsilon}\ell(\gamma_{\varepsilon})\right\vert
_{\varepsilon=0}  &  =\left.  \frac{d}{d\varepsilon}\right\vert _{\varepsilon
=0}%
{\displaystyle\int\nolimits_{0}^{L}}
\left\Vert \gamma_{\varepsilon}^{\prime}(u)\right\Vert du=%
{\displaystyle\int\nolimits_{0}^{L}}
\left.  \frac{d}{d\varepsilon}\right\vert _{\varepsilon=0}\left\Vert
\gamma_{\varepsilon}^{\prime}(u)\right\Vert du\\
&  =%
{\displaystyle\int\nolimits_{0}^{L}}
\left\Vert \gamma^{\prime}\right\Vert ^{-1}V^{\prime}\cdot\gamma^{\prime}du=%
{\displaystyle\int\nolimits_{0}^{L}}
V^{\prime}\cdot\gamma^{\prime}du\\
\left.  \frac{d}{d\varepsilon}\kappa\Gamma_{\varepsilon}(t)\right\vert
_{\varepsilon=0}  &  =\frac{\ell(\gamma)}{L}\cdot\left.  \frac{d}%
{d\varepsilon}\kappa\gamma_{\varepsilon}(t)\right\vert _{\varepsilon
=0}+\left.  \frac{d}{d\varepsilon}\frac{\ell(\gamma_{\varepsilon})}%
{L}\right\vert _{\varepsilon=0}\cdot\kappa\gamma(t)\\
&  =\left(  \gamma^{\prime\prime}\cdot V^{\prime\prime}-2V^{\prime}\cdot
\gamma^{\prime}\right)  (t)+\frac{1}{L}%
{\displaystyle\int\nolimits_{0}^{L}}
V^{\prime}\cdot\gamma^{\prime}du
\end{align*}

\end{proof}

\begin{lemma}
If there exists a $C^{1,1}$ $K$ parametrized by $\gamma:\mathbf{S}%
^{1}\rightarrow K\subset\mathbf{R}^{n}$ satisfying both

(i) $K$ is a relative minimum of $\ell_{e}$ in $\mathbf{R}^{n}$, and

(ii) $NIR(K,\mathbf{R}^{n})=F_{k}(K)=1<\frac{1}{2}DCSD(K),$

then $\exists\delta>0$ such that $\gamma([a,b])$ is a shortest curve in
$\mathcal{C}(\gamma(a),\gamma(b);\gamma^{\prime}(a),\gamma^{\prime}(b))$
whenever $\ell_{ab}(\gamma)\leq\delta.$ This is also true for $K$ with several components.
\end{lemma}

\begin{proof}
By [D2, \emph{Corollary 2}], $\kappa\gamma\equiv1.$ Reparametrize
$\gamma:\mathbf{S}^{1}\cong\mathbf{R}/L\mathbf{Z}\rightarrow K\subset
\mathbf{R}^{n}$, such that $\left\Vert \gamma^{\prime}\right\Vert =\left\Vert
\gamma^{\prime\prime}\right\Vert =1.$ Suppose that such $\delta>0$ does not
exist. $\forall m\in\mathbf{N},\exists a_{m},b_{m}\in\mathbf{S}^{1}$ such that
$0<\left\vert a_{m}-b_{m}\right\vert \leq\frac{1}{m}$ but $\gamma([a_{m}%
,b_{m}])$ is not a shortest curve in $\mathcal{C}(\gamma(a_{m}),\gamma
(b_{m});\gamma^{\prime}(a_{m}),\gamma^{\prime}(b_{m})).$ By [\emph{Proposition
3}, D2], there exists a shortest curve $\theta_{m}$ in $\mathcal{C}%
(\gamma(a_{m}),\gamma(b_{m});\gamma^{\prime}(a_{m}),\gamma^{\prime}(b_{m}))$.
The $C^{1}$ end point data of $\theta_{m}$ and $\gamma([a_{m},b_{m}])$ match.
Let $\gamma_{m}$ be the $C^{1}$ curve obtained from $\gamma$ by removing
$\gamma([a_{m},b_{m}])$ and attaching $\theta_{m}$ in its place. Then,
$\kappa\gamma_{m}\leq1,$ and $L-\frac{1}{m}\leq\ell(\gamma_{m})<L=\ell
(\gamma).$ Hence, it is possible to reparametrize $\gamma_{m}$ uniformly with
a common domain $\mathbf{S}^{1}\cong\mathbf{R}/L\mathbf{Z}$, with $\left\Vert
\gamma_{m}^{\prime}\right\Vert \leq1$ and $\left\Vert \gamma_{m}^{\prime
\prime}\right\Vert \leq1,$ almost everywhere, for sufficiently large $m.$
Hence the convergence $\gamma_{m}\rightarrow\gamma$ can be taken uniformly in
$C^{1}$-sense, by Arzela-Ascoli Theorem and taking a subsequence if it is
necessary. Let $K_{m}=\gamma_{m}(\mathbf{S}^{1}).$ By \emph{Proposition
1(b)},
\[
\lim\inf_{m}\frac{1}{2}DCSD(K_{m})\geq\frac{1}{2}DCSD(K)>F_{k}(K)=1
\]
For sufficiently large $m,$%
\begin{align*}
\frac{1}{2}DCSD(K_{m})  &  >1\\
F_{k}(K_{m})  &  =F_{k}(K)=1\\
NIR(K_{m},\mathbf{R}^{n})  &  =NIR(K,\mathbf{R}^{n})=1\\
\ell(K_{m})  &  <L=\ell(K)\\
\ell_{e}(K_{m})  &  <\ell_{e}(K)
\end{align*}
which contradicts relative minimality of $K$ for $\ell_{e}.$

If $K$ has finitely many components, then at least one of the components of
$K$ contains infinitely many pairs $\{\gamma(a_{m}),\gamma(b_{m})\}$ specified
as above, and the rest of the proof is the same.
\end{proof}

\bigskip

\begin{proof}
(\textbf{THEOREM 1) }Let $n$\ be a dimension such that

(i) every minimizer for the Markov-Dubins problem in $\mathbf{R}^{n}$ is
either a smooth curve with curvature $1$ and positive torsion, or a $C^{1}%
-$concatenation of finitely many circular arcs of curvature $1$ and a line
segment, and

(ii) every $CCC-$curve with the middle arc of length $<\pi$\ is not a minimizer.

First consider the case of a connected $K.$ Suppose that there exists a
relative minimum $K$ of $\ell_{e}$ such that $NIR(K,\mathbf{R}^{n}%
)=F_{k}(K)<\frac{1}{2}DCSD(K).$ Rescale to obtain $F_{k}(K)=1.$ By [D2,
\emph{Corollary 2}], $K$ has constant generalized curvature $\kappa=1.$ By
\emph{Lemma 3}, $\exists\delta>0$ such that $\forall a\in\mathbf{S}^{1},$
$\gamma([a,a+\delta])$ is a shortest curve in $\mathcal{C}(\gamma
(a),\gamma(a+\delta);\gamma^{\prime}(a),\gamma^{\prime}(a+\delta)),$ where
$\gamma:\mathbf{S}^{1}\cong\mathbf{R}/L\mathbf{Z}\rightarrow K\subset
\mathbf{R}^{n}$ is a parametrization with respect to arclength. By the
hypothesis (i), each $\gamma([a,a+\delta])$ is either (a) a smooth curve with
$\kappa=1$ and $\tau>0$ or (b) a $C^{1}-$concatenation of finitely many pieces
each of which is an arc of a circle or a line segment. Since $\kappa=1,$ there
are no line segments. Type (a) curves and type (b) curves do not have any
curve in common, even in part, since one has $\tau>0$ everywhere and the other
is concatenation of planar arcs of circles. $\gamma([a,a+\delta])$ and
$\gamma([a+\frac{\delta}{2},a+\frac{3\delta}{2}])$ have a common piece and
hence they must be of the same type. Inductively, we conclude that either all
of $K$ is a smooth curve with positive torsion or a $C^{1}-$concatenation of
finitely many circular arcs. \emph{Proposition 5} and $F_{k}(K)<\frac{1}%
{2}DCSD(K)$ exclude the smooth case with $\tau>0,$ and imply that $K$ must be
a concatenation of finitely many circular arcs.

We will assume that two successive circular arcs have distinct centers, i.e.
no trivial concatenations. If any of the circular arcs of $K$ has length $\pi$
or more, one can find 2 diametrically opposed points on it, forming a minimal
double critical pair, and $NIR(K,\mathbf{R}^{n})=F_{k}(K)=\frac{1}%
{2}DCSD(K)=1$ contradicting the initial assumption.

This leaves us the final case of $C^{1}-$concatenations with circular arcs of
length $<\pi.$ There must be at least $3$ circular arcs in $K.$ Consider a
parametrization $\gamma:\mathbf{S}^{1}\rightarrow K$ with respect to arclength
such that $\gamma([0,a])$ is a single maximal circular arc of length $a<\pi.$
For $m$ sufficiently large$\mathbf{,}$ $\gamma([-\frac{1}{m},a+\frac{1}{m}])$
is a $CCC$-curve such that the middle arc has length $a<\pi.$ By the
hypothesis (ii), this type $CCC$-sections of $K$ are not minimizers in a
corresponding $\mathcal{C}$. Let $\mathcal{U}$ be an open set in $C^{1}$
topology such that $\gamma\in\mathcal{U}$ and $\ell_{e}(\gamma)\leq\ell
_{e}(\eta),\forall\eta\in\mathcal{U}\cap\lbrack\gamma].$ When one replaces a
non-minimal $CCC$-section with a minimal curve in the same $\mathcal{C}$, then
a priori one can not assume that the new curve is in $\mathcal{U}\cap
\lbrack\gamma]$, and one can not use relative minimality of $K.$ Let
$\theta_{m}$ be any minimizer of $\mathcal{C}(\gamma(-\frac{1}{m}%
),\gamma(a+\frac{1}{m});\gamma^{\prime}(-\frac{1}{m}),\gamma^{\prime}%
(a+\frac{1}{m})).$ Let $\gamma_{m}$ be the $C^{1}$ curve obtained from
$\gamma$ by removing $\gamma([-\frac{1}{m},a+\frac{1}{m}])$ and attaching
$\theta_{m}$ in its place. Then, $\kappa\gamma_{m}\leq1$ and $L-(a+\frac{2}%
{m})\leq\ell(\gamma_{m})<L=\ell(\gamma)\ $where $L-a>0 $ and $L<\infty.$
Hence, for all sufficiently large $m,$ it is possible to reparametrize
$\gamma_{m}$ with a common domain $\mathbf{S}^{1}\cong\mathbf{R}/L\mathbf{Z}$,
such that

1. $\gamma_{m}(-\frac{1}{m})=\gamma(-\frac{1}{m})$

2. $\left\Vert \gamma_{m}^{\prime}\right\Vert =1,$ and $\left\Vert \gamma
_{m}^{\prime\prime}\right\Vert \leq1$ almost everywhere, on $[-\frac{1}%
{m},a+\frac{1}{m}],$ and

3. $\left\Vert \gamma_{m}^{\prime}\right\Vert \leq c_{1}<\infty,$ and
$\left\Vert \gamma_{m}^{\prime\prime}\right\Vert \leq c_{2}<\infty$ almost
everywhere, on $\mathbf{S}^{1}\cong\mathbf{R}/L\mathbf{Z.}$

Clearly, $\left(  \left\Vert \gamma^{\prime}(s)-\gamma^{\prime}(t)\right\Vert
\leq C\left\vert s-t\right\vert ,\forall s,t\in I\right)  \Leftrightarrow$

$\left(  \left\Vert \gamma^{\prime\prime}(s)\right\Vert \leq C\text{ for
almost all }s\in I,\text{ and }\gamma^{\prime}\text{ is absolutely
continuous}\right)  $. By Arzela-Ascoli Theorem, there exists a convergent
subsequence (using the same indices), $\gamma_{m}\rightarrow\gamma_{0}$
converging uniformly on $\mathbf{S}^{1}$ in the $C^{1}$-sense. Then,

1. $\left\Vert \gamma_{0}^{\prime}\right\Vert =1$ on $[0,a],$

2. $\kappa\gamma_{0}\leq1$ on $[0,a],$ since the inequality $\left\Vert
\gamma_{m}^{\prime}(s)-\gamma_{m}^{\prime}(t)\right\Vert \leq\left\vert
s-t\right\vert $ carries to the limit: $\left\Vert \gamma_{0}^{\prime
}(s)-\gamma_{0}^{\prime}(t)\right\Vert \leq\left\vert s-t\right\vert ,$

3. $\gamma_{0}(0)=\gamma(0),$ and hence $\gamma_{0}^{\prime}(0)=\gamma
^{\prime}(0),$ and

4. Since $\gamma_{m}(t_{m})=\gamma(a+\frac{1}{m})$ for some $t_{m}\in
(-\frac{1}{m},a+\frac{1}{m})$, $\exists t_{0}\in\lbrack0,a]$ such that
$\gamma_{0}(t_{0})=\gamma(a)$, and hence $\gamma_{0}^{\prime}(t_{0}%
)=\gamma^{\prime}(a).$

By [\emph{Proposition 1} of D2], the planar circular arc $\gamma([0,a])$ is
the unique minimizer in $\mathcal{C}(\gamma(0),\gamma(a);\gamma^{\prime
}(0),\gamma^{\prime}(a))$ which also contains $\gamma_{0}([0,t_{0}]).$%
\[
a=\ell(\gamma\lbrack0,a])\leq\ell(\gamma_{0}[0,t_{0}])=t_{0}=\lim(t_{m}%
+\frac{1}{m})\leq\lim(a+\frac{2}{m})=a
\]
Consequently, $a=t_{0}$, $\gamma|[0,a]=$ $\gamma_{0}|[0,a]$\ and
$\gamma(\mathbf{S}^{1})=\gamma_{0}(\mathbf{S}^{1})$. Let $K_{m}=\gamma
_{m}(\mathbf{S}^{1}).$ By \emph{Proposition} \emph{1(b)},
\[
\lim\inf_{m}\frac{1}{2}DCSD(K_{m})\geq\frac{1}{2}DCSD(K)>F_{k}(K)=1
\]
For sufficiently large $m,$%
\begin{align*}
\gamma_{m}  &  \in\mathcal{U}\cap\lbrack\gamma_{0}]=\mathcal{U}\cap
\lbrack\gamma]\\
\frac{1}{2}DCSD(K_{m})  &  >1\\
F_{k}(K_{m})  &  =F_{k}(K)=1\\
NIR(K_{m},\mathbf{R}^{n})  &  =NIR(K,\mathbf{R}^{n})=1\\
\ell(K_{m})  &  <L=\ell(K)\\
\ell_{e}(K_{m})  &  <\ell_{e}(K)
\end{align*}
which contradicts relative minimality of $K.$ This shows the nonexistence of
concatenations only with circular arcs of length $<\pi.$ Actually, the
existence of one circular arc of length $<\pi$ actually led to the
contradiction. Since all cases lead to a contradiction, one must have
$NIR(K,\mathbf{R}^{n})=\frac{1}{2}DCSD(K).$

The extension to several component case is straightforward, by
\emph{Proposition 5}, [\emph{Corollary 2} of D2], and the proof of the final
case being a local argument.
\end{proof}

\section{References}

[BS]\qquad G. Buck and J. Simon, \textit{Energy and lengths of knots,
}Lectures at Knots 96, 219-234.

[CKS]\qquad J. Cantarella, R. B. Kusner, and J. M. Sullivan, \textit{On the
minimum ropelength of knots and links, }Inventiones Mathematicae \textbf{150}
(2002) no. 2, p. 257-286.

[CE]\qquad J. Cheeger and D. G. Ebin, \textit{Comparison theorems in
Riemannian geometry, Vol 9, }North-Holland, Amsterdam, 1975.

[Di]\qquad Y. Diao, \textit{The lower bounds of the lengths of thick knots,
}Journal of Knot Theory and Its Ramifications, Vol. \textbf{12}, No. 1 (2003) 1-16.

[Du]\qquad L. E. Dubins, \textit{On curves of minimal length with a constraint
on average curvature and with prescribed initial and terminal positions and
tangents, }Amer. J. Math. \textbf{79} (1957), 497-516.

[D1]\qquad O. C. Durumeric, \textit{Thickness formula and }$C^{1}%
-$\textit{compactness of }$C^{1,1}$\textit{\ Riemannian submanifolds,
}preprint, http://arxiv.org/abs/math.DG/0204050

[D2]\qquad O. C. Durumeric, \textit{Local structure of the ideal shapes of
knots, preprint, }revised version, December 2005\textit{.}

[FHW]\qquad M. Freedman, Z.-X. He, and Z. Wang, \textit{Mobius energy of knots
and unknots}, Annals of Math, \textbf{139} (1994) 1-50.

[GL]\qquad O. Gonzales and R. de La Llave, \textit{Existence of ideal knots,}
J. Knot Theory Ramifications, \textbf{12} (2003) 123-133.

[GM]\qquad O. Gonzales and H. Maddocks, \textit{Global curvature, thickness
and the ideal shapes of knots}, Proceedings of National Academy of Sciences,
\textbf{96 }(1999) 4769-4773.

[GMSM]\qquad O. Gonzales, H. Maddocks, F. Schuricht, and H. von der Mosel,
Global curvature and self-contact of nonlinearly elastic curves and rods,
\textit{Calc. Var. Partial Differential Equations}, \textbf{14} (2002) 29-68.

[Ka]\qquad V. Katrich, J. Bendar, D. Michoud, R.G. Scharein, J. Dubochet and
A. Stasiak, \textit{Geometry and physics of knots}, Nature, \textbf{384}
(1996) 142-145.

[L]\qquad A. Litherland, \textit{Unbearable thickness of knots}, preprint.

[LSDR]\qquad A. Litherland, J Simon, O. Durumeric and E. Rawdon,
\textit{Thickness of knots}, Topology and its Applications, \textbf{91}(1999) 233-244.

[N]\qquad\ A. Nabutovsky, \textit{Non-recursive functions, knots
\textquotedblleft with thick ropes\textquotedblright\ and self-clenching
\textquotedblleft thick\textquotedblright\ hyperspheres}, Communications on
Pure and Applied Mathematics\textit{, }\textbf{48} (1995) 381-428.

[M]\qquad A. A. Markov, \textit{Some examples of the solution of a special
kind of problem in greatest and least quantities, }(in Russian), Soobbshch.
Karkovsk. Mat. Obshch. \textbf{1} (1887), 250-276.

[RS]\qquad E. Rawdon and J. Simon, \textit{Mobius energy of thick knots},
Topology and its Applications,\textbf{125} (2002) 97-109.

[ReSh]\qquad J. A. Reeds and L. A. Shepp, \textit{Optimal paths for a car that
goes both forwards and backwards, }Pacific J. Math. \textbf{145} (1990), 367-393.

[S]\qquad H. J. Sussmann, \textit{Shortest 3-dimensional paths with a
prescribed curvature bound,} Proceedings of the 34th Conference on Decision \&
Control, New Orleans, LA-December 1995, IEEE Publications, New York, 1995, 3306-3312.

\end{document}